\documentclass[10pt,a4paper]{article}
\usepackage[T1]{fontenc}
\usepackage{amsfonts,amsmath,amsthm,amssymb,mathrsfs}
\usepackage{graphicx,geometry}
\usepackage{array}
\usepackage[english]{babel}
\usepackage{enumerate}

\usepackage{algorithm}
\usepackage{algpseudocode}
\usepackage{url}

\algrenewcommand\algorithmicwhile{\textbf{While}}
\algrenewcommand\algorithmicfor{\textbf{For}}
\algrenewcommand\algorithmicdo{\textbf{Do}}
\algrenewcommand\algorithmicif{\textbf{If}}
\algrenewcommand\algorithmicthen{\textbf{Then}}
\algrenewcommand\algorithmicelse{\textbf{Else}}
\algrenewcommand\algorithmicend{\textbf{End}}
\algrenewcommand\algorithmicreturn{\textbf{Return}}

\newcommand{\RR}{\mathbb{R}}
\newcommand{\QQ}{\mathbb{Q}}
\newcommand{\CC}{\mathbb{C}}
\newcommand{\PP}{\mathbb{P}}

\newcommand{\Id}{\mathrm{Id}}

\newtheorem{theorem}{Theorem}[section]
\newtheorem{lemma}[theorem]{Lemma} 
\newtheorem{proposition}[theorem]{Proposition}

\newtheorem{remark}[theorem]{Remark}  
\newtheorem{example}[theorem]{Example} 

\newcommand{\pcoor}[1]{%
	\begingroup\lccode`~=`: \lowercase{\endgroup
		\edef~}{\mathbin{\mathchar\the\mathcode`:}\nobreak}%
	(% opening symbol
	\begingroup
	\mathcode`:=\string"8000
	#1%
	\endgroup 
	)% closing symbol
}

\title{
	Using Algebraic Geometry to Reconstruct a Darboux Cyclide from a Calibrated Camera Picture
}
\author{Eriola Hoxhaj$^1$, Jean Michel Menjanahary$^2$, Josef Schicho$^3$}
\date{$^1$ Johannes Kepler University/ RISC, Austria\\[1ex]
	$^2$ Vilnius University/MIF/ Institute of Informatics, Lithuania\\[1ex]
	$^3$ Johannes Kepler University/ RISC, Austria} % <--- use \date for affiliation
\begin{document}
	\maketitle
	
	\begin{abstract}
		The task of recognizing an algebraic surface from a single apparent contour can be reduced to
		the recovering of a homogeneous equation in four variables from its discriminant.
		In this paper, we use the fact that Darboux cyclides have a singularity along the absolute conic in order to recognize them up to Euclidean similarity transformations.
	\end{abstract}
	
	\section{Introduction}
	\label{intro}
	Can we obtain complete spatial information about an algebraic surface in 3-space from a 2D picture?
	It seems counter-intuitive at first.
	If the surface is smooth, then it can be reconstructed from a single apparent contour up to a projective equivalence that fixes all lines passing through the camera location~\cite{Almeida,Forsyth}.
	D'Almeida's algorithm has been generalized in \cite{GLSV} to surfaces with ordinary singularities (nodal curves, transversal triple points, and pinch points).
	
	In this paper, we focus on the reconstruction from a single view by a calibrated camera of a Darboux cyclide. These are surfaces in 3-space that have at least two circles passing to a generic point on it (see \cite{Niels,Lubbes}). In \cite{Fou:04}, special Darboux cyclides called Dupin cyclides are used to blend between quadratic surfaces.
	
\begin{figure}
		\includegraphics[width=0.4\textwidth]{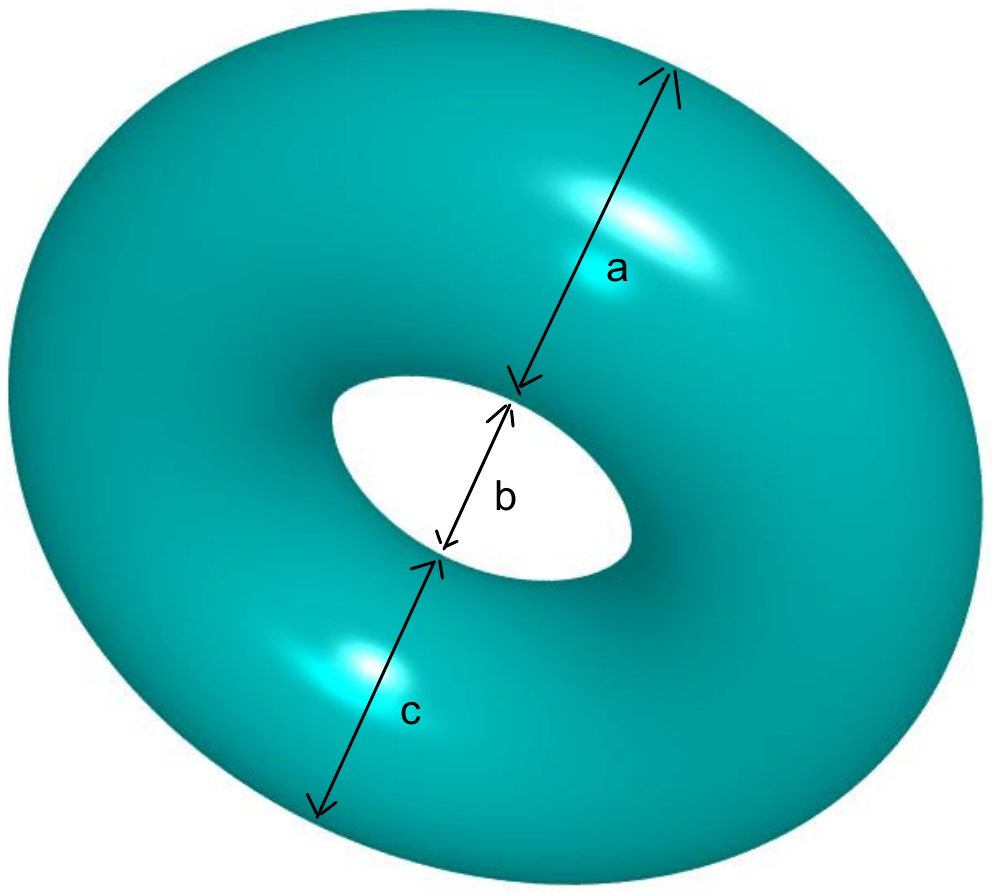} 
		\includegraphics[width=0.4\textwidth]{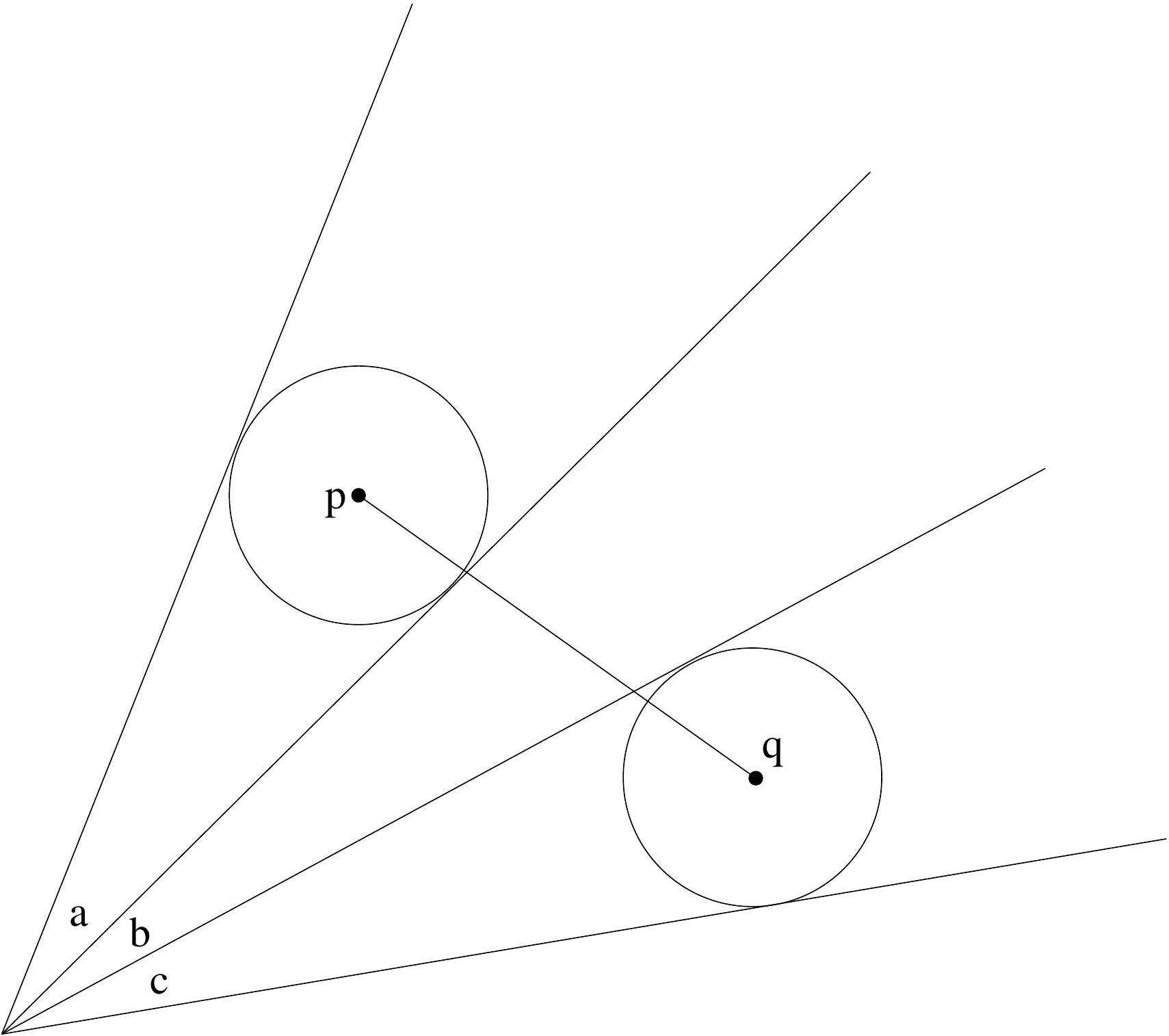}
		\caption{Reconstruction of a torus from a single picture. From the photo with a calibrated camera, the three angles $a,b,c$ are
			determined and drawn in a plane. Next, we draw a circle tangential to the rays forming the angle $a$. Next, we draw a second circle 
			with equal radius tangential to the rays forming the angle $c$. Next, we rotate both circles around the symmetry line
			of the two midpoints $p,q$ and obtain the torus.}
		\label{fig:torus}
	\end{figure}

	A picture from a single view by a calibrated camera allows to determine viewing angles (see \cite{Hartley_Zisserman}). It is easy to reconstruct the equation of a sphere from a single view: the viewing angle already determines everything. 
	For a ring torus, the situation is slightly more complicated but still doable -- see Figure~\ref{fig:torus} for a geometric construction.
	Here, we give an algorithm that reconstructs a Darboux cyclide from the apparent contour of a single view with a calibrated camera,
	assuming that the camera is in a general position.
	Up to scaling, there are only finitely many solutions.
	In particular, if we have one Darboux cyclide with given apparent contour, then we can construct a second one by inversion at a sphere centered at the camera position.
	
	Algorithmically, our results are strongly based on the ideas in \cite{Almeida}, and the theorems and algorithms in \cite{GLSV}.
	The stepping stone is the contour in projective 3-space, including the non-reduced structure in the Euclidean absolute conic.
	To compute it, we collect local contributions from special points in the apparent contour: 
	non-smooth points of the apparent contour that do not arise by projecting special singular points of the Darboux cyclide.
	The local contributions are generators for the local conductor ideal.
	In many cases, the formula for the conductor ideal is well-known; there is one case which is new (Lemma~\ref{lem:new}).
	The non-uniqueness of the solution is explained by the following uncertainty for a non-smooth point $q$ of the apparent contour.
	In order to proceed with the computation, we need to know -- or to assume -- whether $q$ is special or whether it is the image of a singular point of the surface. 
	Depending on the assumptions we make at this step, we get different solutions, or no solution in case the assumptions are inconsistent.
	
	The contents of the paper are organized as following.  In Section~\ref{sec:1}, we give a mathematical definition of a pinhole camera 
	and the calibration we consider on a pinhole camera. Then we formulate accordingly the problem of reconstruction of an 
	algebraic surface in space from a single view. 
	% All possible solutions of this reconstruction is discussed as well as the 
	% sub-solutions obtained by using the additional information carried by the camera calibration. 
	In Section~\ref{sec:2}, we recall the definition of Darboux cyclides and their relationship with the absolute conic.
	We investigate the 
	singularities of Darboux cyclides; this is needed later for the reconstruction procedure.  
	In Section~\ref{sec:3}, we give an analysis of the contour and apparent contour.
	We especially discuss the singularities of both curves.
	%The contour is known as the ``ramification curve'' and the apparent contour is known as the ``branching curve'', following standard terminology in algebraic geometry.
	%In Section 3, we explore the input and the first output of our algorithm reconstruction: 
	%the input for Darboux cyclides is a planar curve of degree 6 or 8 with an additional non-real component curve 
	%that is an invisible 2d image of the special conic. The output will be a space curve that 
	%is projected to the given input curve. 
	Section~\ref{sec:4} contains our algorithm for the construction of the contour from the given apparent contour.
	Especially, this section contains formulas for the local contributions to the conductor mentioned above.
	%The ideas and proofs of the reconstruction of the space curve is explored 
	%in Section 4, which is based on well-known results and previous works in this direction: we collect the local 
	%contribution of the singularities on the input curve that are obtained by the projection and apply the theories in 
	%order to reconstruct the space curve. 
	Section~\ref{sec:5} explains the last step of the algorithm, namely the construction of the equation of the surface
	when the equation of the contour is already known. 
	
	\section{The Calibrated Camera}
	\label{sec:1}
	
	The mathematical model of a {\em pinhole camera} is a linear projective map $\pi:(\PP^3\setminus\{o\})\to\PP^2$ given by a $3\times 4$ matrix,
	where $o\in\PP^3$ is the position of the camera; its representation by homogeneous coordinates spans the kernel of the matrix.
	In this paper, we study the recognition of an algebraic surface $\mathrm{S}:F(x,y,z,w)=0$ by a single view. So, we assume without loss of generality
	that $o=(0:0:0:1)$ and that $\pi$ is the projection map $(x:y:z:w)\mapsto (x:y:z)$. 
	The {\em contour} of $\mathrm{S}$ is the space curve defined by $F=\partial_w F=0$ -- geometrically, this is the locus of all points
	with a tangent passing through $o$. 
	The {\em apparent contour} of $\mathrm{S}$ is the image of the contour. It is an algebraic curve defined by the 
	discriminant $\Delta_w(F)$. Single view surface recognition is the task of reconstructing $F$ from its discriminant $\Delta_w(F)$.
	For any $a,b,c\in\RR$ and $d\in\RR^\ast$,
	the polynomials $F(x,y,z,w)$ and $F(x,y,z,ax+by+cz+dw)$ have the same discriminant up to a constant factor. 
	Therefore, single view recognition is at best only possible up to the group of projective transformations of the form
	$(x:y:z:w)\mapsto (x:y:z:ax+by+cz+dw)$. These are the projective transformations fixing all lines passing through $o$ (and hence $o$ as well).
	
	The camera is calibrated if and only if the image plane is an {\em elliptic plane}. This means that it comes with a metric, called
	{\em elliptic distance}. For any two points $a,b\in\PP^3\setminus\{p\}\to\PP^2$, the elliptic distance of $\pi(a)$ and $\pi(b)$ 
	is equal to the angle $\angle(apb)$. The metric is determined by the {\em absolute conic} $\mathrm{C}_E$ of the elliptic plane, an irreducible
	conic without real points. Here, we will assume that $\mathrm{C}_E$ has equation $x^2+y^2+z^2=0$. The elliptic distance of two points  $\pi(a)$ and
	$\pi(b)$ can then be computed as half of the absolute value of the arctangent of $\frac{i(x^2-1)}{x^2+1}$, where $x$ is the cross ratio (see \cite[p.45]{Hartley_Zisserman}) of $\pi(a)$, $\pi(b)$, and the other two intersection points of the line $\pi(a)\pi(b)$ with $\mathrm{C}_E$, and $i$ is the imaginary unit.
	
	\begin{remark} \rm A well-known way to ``implement'' a calibrated camera
		is to mark on the footpoint of $o$ to the image plane on the picture and to specify the elliptic distance of this footpoint to a single other
		point. The footpoint has elliptic distance $\pi/2$ to all points on the line of infinity, and one can show that the position of the
		footpoint, the position of the line at infinity, and a single distance between the footpoint and any finite point uniquely determines
		the absolute conic of the elliptic plane. 
	\end{remark}
	
	The consequence of our assumption that the camera given by the map $(x:y:z:w)\mapsto (x:y:z)$ is calibrated is the following: 
	the 3-space is not just a projective 3-space, but it also has the structure of a Euclidean space, i.e. there is a Euclidean metric
	defined on the set of finite points. In terms of projective geometry, the Euclidean metric is defined up to scaling by the absolute
	conic $\mathrm{C}_A$, a conic on the plane of the infinite plane without real points. In this paper, $\mathrm{C}_A$ is the conic $x^2+y^2+z^2=w=0$.
	In general, it is required that $\pi$ maps the Euclidean absolute conic to the absolute conic of the elliptic plane; and this is
	obviously true in our setup.
	
	A calibrated camera should, in principle, give more precise information than a camera which is not calibrated. For instance, we may
	hope for an answer that is unique not just up to projective transformations of the form $(x:y:z:w)\mapsto (x:y:z:ax+by+cz+dw)$,
	but more specifically unique up to scaling $(x:y:z:w)\mapsto (x:y:z:dw)$ (recall $d\ne 0$). The reason for such a hope is that
	all maps of the form $(x:y:z:w)\mapsto (x:y:z:ax+by+cz+dw)$ that preserve the Euclidean absolute conic also preserve the plane
	at infinity, and they must therefore be scalings.
	
	The hope expressed above has a chance to be fulfilled only if the surface to be recognized has some relation to the Euclidean absolute conic.

	\section{Darboux Cyclides}
	\label{sec:2}
	
	A \emph{Darboux cyclide} is an algebraic surface of degree four, obtained by intersecting the 3-dimensional M\"obius sphere $x^2+y^2+z^2+t^2=w^2$ with a quadratic surface in $\PP^4$. To obtain a model in $\PP^3$, we project stereographically from $(0:0:0:1:1)$. The three-dimensional image is defined by a polynomial of the form
	\begin{align}
		F(x,y,z,w) = A^2 + 2ALw + Q w^2,
	\end{align}
	where $A:=x^2+y^2+z^2$, $L$ is linear in $x,y,z$ and $Q$ is homogeneous quadratic polynomial in $x,y,z,w$. 
	To avoid degenerate cases, we will assume that the polynomial $F$ is absolutely irreducible.
	This excludes cubic and quadratic Darboux cyclides (see \cite{Pottmann:12,Zhao:19} for a discussion of these degenerate cases) 
	that would have the infinite plane $w=0$ as a second component.
	
	Since $F$ is contained in the ideal $\langle A,w\rangle_{\QQ[x,y,z,w]}^2$, the absolute conic $\mathrm{C}_A$ 
	-- which is defined by $A=w=0$ -- is at least a double curve of the Darboux cyclide. 
	We need to make this statement a little bit more specific. 
	A {\em nodal curve} of a surface $\mathrm{S}$ is an irreducible curve $\mathrm{C}$ in the singular locus of $F$ such that the intersection of $\mathrm{S}$ with a generic transversal plane $\mathrm{P}$ has an ordinary double point (also known as ``node'') at any intersection of $\mathrm{P}$ and $\mathrm{C}$.
	A {\em cuspidal curve} of a surface $\mathrm{S}$ is an irreducible curve $\mathrm{C}$ in the singular locus of $F$ such that the intersection of $\mathrm{S}$ with a generic plane $\mathrm{P}$ has a cusp at any intersection of $\mathrm{P}$ and $\mathrm{C}$.
	
	\begin{proposition}
		Let $\mathrm{D}\subset\PP^3$ be an irreducible Darboux cyclide. 
		Then the absolute conic $\mathrm{C}_A$ is either a nodal or a cuspidal curve of $\mathrm{D}$. 
		Also, the absolute conic is the only singular curve of $\mathrm{D}$.
	\end{proposition}
	
	\begin{proof}
		A generic plane $\mathrm{P}$ intersects $\mathrm{D}$ in a plane quartic curve $\mathrm{Q}$. 
		By Bertini's Theorem \cite[p.137]{griffiths:14}, the quartic $\mathrm{Q}$ is irreducible. 
		Since $\mathrm{C}_A$ has no real points, the singularities of $\mathrm{Q}$ come in pairs of complex conjugates. 
		By the genus formula for plane algebraic curves (see \cite{Hartshorne:77}, exercise IV.1.8a), 
		and by the fact that the genus cannot be negative, the sum of all $\delta$-invariants cannot be bigger than $3$. 
		On the other hand, conjugate singularities have the same $\delta$-invariant. 
		Hence, we have exactly one pair of conjugated singularities, which have $\delta$-invariant equal to one. 
		The only singularities with $\delta$-invariant equal to one are nodes and cusps.
		This shows that $\mathrm{C}_A$ is nodal or cuspidal.
		
		Assume, indirectly, that there exists another singular curve $\mathrm{B}$. 
		It would intersect $\mathrm{P}$ in $\deg(\mathrm{B})$ many points, and these points would be singular points of $\mathrm{Q}$.
		Since $\mathrm{Q}$ has at most three singular points, it follows that $\deg(\mathrm{B})=1$, i.e., $\mathrm{B}$ is a line.
		The line $\mathrm{B}$ intersects the infinite plane in a single point $p$, which must be real because otherwise the intersection
		would also contain the conjugate point.
		On the other hand, the infinite plane intersects $\mathrm{D}$ only in the absolute conic, with multiplicity 2.
		This is a contradiction because the absolute conic has no real points. 
	\end{proof}
	
	\begin{proposition} \label{prop:no3}
		An irreducible Darboux cyclide does not have triple or quadruple points. 
	\end{proposition} 
	
	\begin{proof}
		If a Darboux cyclide $\mathrm{D}$ would have a quadruple point, then it would be a cone, which is obviously not the case.
		
		Assume, indirectly, that $\mathrm{D}$ has a triple point $p$. For any point $q\in \mathrm{C}_A$, the line $\mathrm{L}_{pq}$ must be contained in $\mathrm{D}$. 
		Then $\mathrm{D}$ would contain a whole quadric cone, or plane in case $p$ is lying on the infinite plane. This contradicts irreducibility. 
	\end{proof}

	\section{The contour and apparent contour}
	\label{sec:3}
	
	Assume that $\mathrm{D}$ is an irreducible Darboux cyclide given by the equation $F(x,y,z,w)=0$.
	We consider the image of $\mathrm{D}$ by a generic projection. 
	Still, we want to assume that the projection is $\pi:\PP^3\dashrightarrow\PP^2$, $(x:y:z:w)\mapsto(x:y:z)$. 
	In order to achieve this, we translate the center of projection, which is a generic point in $\PP^3$, to the point $o=(0:0:0:1)$.
	The equation after the translation is again called $F$, admittedly a slight abuse of notation.
	As it has been pointed out by a reviewer, this translation is only possible if we do have a central projection and not a parallel projection,
	with center at infinity. In this paper, we do not deal with the problem of reconstructing a Darboux cyclide from an image under parallel projection.
	
	The {\em contour} $\mathrm{R}$  is the common zero set of $F$ and $\partial_w F$.
	Singular curves of $\mathrm{D}$ are always components of $\mathrm{R}$, with some multiplicity.
	The notion of multiplicity for a space curve requires an explanation: 
	we consider $\mathrm{R}$ as the 1-cycle defined as the intersection of $\mathrm{D}$ and the cubic surface defined by $\partial_w F=0$.
	As such, it is a formal sum of irreducible curves with some multiplicity.
	The sum of degree times multiplicity over all components is equal to 12, by Bezout's Theorem.
	In our situation, $\mathrm{C}_A$ is a component of $\mathrm{R}$ of degree~2.
	Its multiplicity is $2$ if $\mathrm{C}_A$ is nodal and $3$ if $\mathrm{C}_A$ is cuspidal.
	Of course, there must be another component.
	
	\begin{proposition}
		Besides $\mathrm{C}_A$, the contour has exactly one other component.
		Its multiplicity is one.
		Its degree is 8 if $\mathrm{C}_A$ is nodal and 6 if $\mathrm{C}_A$ is cuspidal.
	\end{proposition}
	
	\begin{proof}
		The curve is a generic element in the linear system cut out by the partial derivatives of $F$. 
		By Bertini's Theorem, it is smooth outside the base locus, which is the singular locus of $\mathrm{D}$. 
		It follows that $\mathrm{C}_A$ is the only multiple component.
		
		The {\em polar map}, defined by the partial derivatives, maps $\mathrm{D}$ to its dual (see \cite[Section~1.2]{cag}).
		Since $\mathrm{D}$ is not a ruled surface, the dual is two-dimensional.
		By Bertini's Theorem again, it follows that the preimage of a generic plane section is irreducible.
		
		The degree of the second component can be computed by Bezout's Theorem. 
		It is $12-2\cdot 2=8$ in the nodal case, and $12-2\cdot 3=6$ in the cuspidal case. 
	\end{proof}
	
	The {\em apparent contour} $\mathrm{B}\subset\PP^2$ is defined as the image $\pi(\mathrm{R})$ of the contour.
	The elliptic absolute conic $\mathrm{C}_E$ with equation $A=0$ is a double (in the nodal case) or a triple (in the cuspidal case) component.
	The second component has degree $8$ or $6$. We denote it by $\mathrm{B}_1$.
	This is the curve which is visible in the picture.
	The second component of $\mathrm{R}$ -- the one that maps to $\mathrm{B}_1$ -- is denoted by $\mathrm{R}_1$.
	
	The equation of $\mathrm{B}$ is the input of our reconstruction algorithm. 
	It is a polynomial $U\in\QQ[x,y,z]$ of degree $12$ that factors into $U=U_1A^2$ or $U=U_1A^3$.
	The first step of our algorithm is to reconstruct the contour $\mathrm{R}$.
	It will be necessary to reconstruct $\mathrm{R}$ as a scheme, not just as a 1-cycle: we need the ideal of $\mathrm{R}$. 
	Note that we already know that $\mathrm{R}$ is generated by a quartic and a cubic form in $\QQ[x,y,z,w]$.
	
	For any $q\in \mathrm{B}$, we say that the map $\pi|_\mathrm{R}:\mathrm{R}\to \mathrm{B}$ is locally an isomorphism over $q$ if 
	there is a neighborhood $U$ of $q$ such that the map maps the inverse image of $U$ isomorphically to $U$.
	This implies in particular that the preimage of $q$ has only a single point (but maybe with non-reduced scheme structure).
	
	Also, we call a non-isolated singular point $q$ on a nodal curve {\em special} iff the singular point
	obtained by intersection with a generic plane through $q$ is not a node (but more complicated). Similarly, we call a non-isolated singular point $q$ on a cuspidal curve {\em special} iff the singular point
	obtained by intersection with a generic plane through $q$ is not a cusp.
	
	\begin{example}\rm
		The line $x=y=0$ is a nodal curve on the ``Whitney umbrella'' with equation $x^2-y^2z=0$. 
		The point $q=(0,0,0)$ is special: a generic plane section through $q$ has a cusp at $q$.
		A special point with this property is also called a {\em pinch point}.
	\end{example}
	
	\begin{proposition} \label{prop:iso}
		The map $\pi|_\mathrm{R}$ is locally an isomorphism over the following points in $\mathrm{B}$:
		\begin{description}
			\item[\rm a)] smooth points of $\mathrm{B}_1$ outside of $\mathrm{C}_E$;
			\item[\rm b)] smooth points of $\mathrm{C}_E$ outside of $\mathrm{B}_1$;
			\item[\rm c)] images of isolated singular points of $\mathrm{D}$;
			\item[\rm d)] images of special non-isolated singular points of $\mathrm{D}$.
		\end{description}
	\end{proposition}
	
	\begin{proof}
		a) The map $\pi|_{\mathrm{R}_1}$ is generically injective, and it is locally injective for all smooth points in $\mathrm{B}_1$.
		Outside of $\mathrm{C}_A$ it coincides with $\pi|_\mathrm{R}$.
		
		c,d) Let $q\in \mathrm{R}$ be an isolated singularity of $\mathrm{D}$ or a special non-isolated singularity. 
		By Proposition~\ref{prop:no3}, $q$ is a double point of $\mathrm{D}$. 
		Since the projection is generic, the fiber of the image point, i.e., the line $qo$, intersects $\mathrm{D}$ in three points, 
		namely $q$ itself and two nonsingular points that are not in $\mathrm{R}$.
		We introduce an analytic coordinate system $(x,y,z)$ around $p$ such that the projection is $(x,y,z)\mapsto (x,y)$ and $q=(0,0,0)$.
		In this coordinate system, the equation of $\mathrm{D}$ is a quadratic Weierstrass polynomial in $z$, that is, a quadratic polynomial
		$\bar{F}(x,y,z)=z^2+F_1(x,y)z+F_2(x,y)$ with $F_1,F_2\in\CC[[x,y]]$
		(Note that the coefficients might be complex because it could be that $q$ is not a real point).
		The contour $\mathrm{R}$ is locally defined by the equations $\bar{F}=\partial_z\bar{F}=0$, which is equivalent to
		$z^2+F_1z+F_2=2z+F_1=0$. 
		The apparent contour $\mathrm{B}$ is locally defined by the discriminant $F_1^2-4F_2=0$.
		So, the projection restricted to $\mathrm{B}$ has local inverse $(x,y)\mapsto\left(x,y,-\frac{F_1(x,y)}{2}\right)$.
		
		b) The proof of (c,d) is also valid for almost all points of $\mathrm{C}_E$, namely for those points $q$ such that the line $qo$ intersects $\mathrm{D}$ in three points.
		If the line intersects $\mathrm{D}$ in only two or one point, then the order of the discriminant is larger than 2 (in the nodal case) respectively 3 (in the cuspidal case).
		Since $\mathrm{C}_E$ is smooth, the order can only be higher if we have a point of intersection with $\mathrm{B}_1$. 
	\end{proof}
	
	Proposition~\ref{prop:iso} indicates that the unknown curve $\mathrm{R}$ is not so different from the known curve $\mathrm{B}$.
	Now we describe the points more closely for which a difference occurs.
	
	\begin{proposition} \label{prop:noniso}
		Assume that $q\in \mathrm{B}$ is a point over which the map $\pi|_\mathrm{R}$ is not locally an isomorphism.
		Then we have one of the following cases.
		\begin{description}
			\item[\rm a)] $q$ is a node on $\mathrm{B}_1$ and does not lie on $\mathrm{C}_E$. 
			The fiber $\pi^{-1}(q)$ intersects $\mathrm{R}_1$ in two smooth points.
			\item[\rm b)] $q$ is a cusp on $\mathrm{B}_1$ and does not lie on $\mathrm{C}_E$.
			The fiber $\pi^{-1}(q)$ intersects $\mathrm{R}$ in a unique smooth point with multiplicity~2.
			\item[\rm c)] $q$ is a transversal intersection of $\mathrm{B}_1$ and $\mathrm{C}_E$.
			The fiber $\pi^{-1}(q)$ intersects $\mathrm{R}_1$ in a smooth point and $\mathrm{C}_A$ in a different point.
			The second point is not a special singularity.
			\item[\rm d)] 
			$q$ is a tangential intersection of $\mathrm{B}_1$ and $\mathrm{C}_E$, with intersection multiplicity $2$ or $3$ depending
			whether $\mathrm{C}_A$ is nodal or cuspidal.
			The fiber $\pi^{-1}(q)$ intersects $\mathrm{R}$ in a single point in the intersection of $\mathrm{R}_1$ and $\mathrm{C}_A$.
			This point is not a special singularity and not a singular point of $\mathrm{R}_1$.
		\end{description}
	\end{proposition}
	
	\begin{proof}
		First, assume that $q$ does not lie on $\mathrm{C}_E$. 
		By Proposition~\ref{prop:iso}, the fiber $\pi^{-1}(q)$ intersects $\mathrm{R}_1$ only in smooth points of $\mathrm{D}$. 
		The apparent contour of a smooth surface under generic projections has only nodes and cusps (see \cite{Ciliberto_Flamini}). This
		shows that we have (a) or (b).
		
		Second, assume that $q$ does lie on $\mathrm{C}_E$.
		The center of projection $o$ does not lie on the plane at infinity.
		Therefore, the fiber $\pi^{-1}(q)$ (a line through $o$) intersects $\mathrm{C}_A$ only once, and the intersection is transversal.
		By Proposition~\ref{prop:iso}, the preimage in $\mathrm{C}_A$ is not a special singularity of $\mathrm{D}$.
		By Proposition~\ref{prop:iso}(b), the fiber also intersects $\mathrm{R}_1$.
		By Lemma~2.3 and Lemma~2.4 in \cite{GLSV}, there is only a single intersection with $\mathrm{R}_1$, and it is transversal.
		We distinguish three subcases.
		
		Subcase 1: the preimages in $\mathrm{R}_1$ and in $\mathrm{C}_A$ are distinct. By Lemma~2.5 in \cite{GLSV}, the tangents at the two preimages to $\mathrm{R}_1$ and to $\mathrm{C}_A$
		are not coplanar. We obtain (c).
		
		Subcase 2: the fiber meets $\mathrm{R}_1$ and $\mathrm{C}_A$ in the same point, and $\mathrm{C}_A$ is nodal.
		By \cite[Proposition 2.1]{GLSV}, the point $q$ is a smooth point in both curves $\mathrm{C}_E$ and $\mathrm{B}_1$, and the two curves meet with
		intersection multiplicity 2.
		
		Subcase 3: the fiber meets $\mathrm{R}_1$ and $\mathrm{C}_A$ in the same point $p$, and $\mathrm{C}_A$ is cuspidal. Locally around $p$, the surface $\mathrm{D}$ is
		a cylinder over the cuspidal curve with equation $y^2-z^3=0$. 
		Since $p$ is also in $\mathrm{R}_1$, the center of projection lies in the plane of maximal contact $y=0$.
		This holds only for $p$, not for other points on $\mathrm{C}_A$ close to $p$. If $p_x$ parametrizes $\mathrm{C}_A$ so that $p_0=p$, then the 
		plane of maximal contact changes with $y$ and passes through $o$ exactly at $x=0$. The situation is described in a different
		local coordinate system as follows:
		the projection map is $(x,y,z)\mapsto (x,y)$, and the surface $\mathrm{D}$ has equation $z^3-(y-xz)^2=0$. The cuspidal curve $\mathrm{C}_A$
		has equation $x=z=0$, and it projects to the curve $\mathrm{C}_E$ with equation $y=0$, which appears as a triple component of the discriminant. The 
		curve $\mathrm{R}_1$ has parameterization $(x,y,z)=\left(x,\frac{4x^3}{27},\frac{4x^2}{9}\right)$, and its image $\mathrm{B}_1$ has equation
		$4x^3-27y=0$ and intersects $\mathrm{C}_E$ with intersection multiplicity three. 
	\end{proof}
	
	\begin{remark} \label{rem:ambi} \rm
		It is not always possible to infer the type of singularity of the contour from the type of singularity of the apparent contour. 
		For instance, a transversal intersection of $\mathrm{C}_E$ and $\mathrm{B}_1$ could be the image of two distinct points as in Proposition~\ref{prop:noniso} case (c), or the image of a pinch point.
		
		As a consequence, the recognition does not always have a unique solution. We can say that the number of solutions, up to scalings, is finite, since for each singularity of the contour we only have a binary choice. In particular, the inversion of a Darboux cyclide at a sphere with center $o$ has the same apparent contour.
		There are several suggestions how to obtain additional information in order to have only a unique answer:
		we can assume that the solution is defined over $\QQ$ and therefore the finite choice which we have to make must be invariant under conjugation by the Galois group of Galois closure of the residue field at the singularity of the contour, or we can use a second view from a different camera position. 
		It would be nice to have a theoretical result giving necessary criteria for uniqueness; unfortunately, we do not have such a result.
	\end{remark}
	
	\section{The conductor of the apparent contour in the contour}
	\label{sec:4}
	
	Let $q$ be a point in $\mathrm{B}$. 
	The local ring $E_q$ of $\mathrm{B}$ at $q$ is defined as the ring of all regular functions on $\mathrm{B}$ defined in some neighborhood of $q$.
	Algebraically, it is the quotient of a two-dimensional regular local ring by a local equation of the discriminant.
	Let $q_1,\dots,q_k\in \mathrm{R}$ be the points in the fiber (by now we know that $k=1$ or $k=2$). 
	We define $F_q$ as the ring of all regular functions on $\mathrm{R}$ defined in neighborhoods of $q_1,\dots,q_k$. 
	If $k=1$, then $F_q$ is the local ring of $\mathrm{R}$ at $q_1$; otherwise, $F_q$ is only semi-local.
	
	The projection $\pi$ induces an injective ring homomorphism $E_q\hookrightarrow F_q$.
	Note that the ring homomorphism is finite: over $F_q$, the ring $E_q$ is generated by a single element -- the vertical coordinate -- 
	that fulfills an integral equation, namely a local equation for $\mathrm{D}$.
	Recall that the total fraction ring of a ring is the localization at all elements that are not zero divisors.
	%In particular, the total fraction ring of the local ring of a scheme is the set of regular functions defined somewhere on every irreducible component.
	Because the projection map is an isomorphism over almost all points, the ring inclusion $E_q\hookrightarrow F_q$ extends to an isomorphism of total fraction rings.
	We denote the total fraction ring of $E_q$ by $K_q$; the total fraction ring of $F_q$ is also identified with $K_q$.
	
	The {\em conductor} $C_q$ of $E_q$ in $F_q$ is defined as the set of all $a\in E_q$ such that $aF_q\subseteq E_q$.
	It is an ideal in both rings $E_q$ and $F_q$.
	Our interest in the conductor comes from the fact that the conductor is isomorphic to $F_q$ as an $E_q$-module.
	Here is the precise statement; the proof can be found in \cite{GLSV}.
	
	\begin{lemma} \label{lem:sheaves}
		Let $\mathrm{X},\mathrm{Y}$ be complete intersections and let $f:\mathrm{X}\to \mathrm{Y}$ be a finite map which is an isomorphism over all generic points of $\mathrm{Y}$.
		Let $\omega^0_\mathrm{X},\omega^0_\mathrm{Y}$ be the dualizing sheafs of $\mathrm{X},\mathrm{Y}$.
		Let ${\cal C}\subset{\cal O}_\mathrm{X}$ be the sheaf of conductors of ${\cal O}_\mathrm{Y}$ in $f^\ast{\cal O}_\mathrm{X}$.
		Then $f_\ast\omega^0_\mathrm{X}\cong {\cal C}{\otimes}\omega^0_\mathrm{Y}$.
	\end{lemma}
	
	The rational map $\mathrm{B}\dashrightarrow \mathrm{R}$ which is inverse to $\pi|_\mathrm{R}$ is defined by the global sections of $(\pi|_\mathrm{R})_\ast{\cal O}_\mathrm{R}(1)$.
	For a complete intersection in $\PP^n$, the dualizing sheaf is the twisted sheaf with degree shift equal to the sum of the
	degrees of the intersecting hypersurfaces minus $n-1$.
	By Lemma~\ref{lem:sheaves}, we can replace the sheaf $(\pi|_\mathrm{R})_\ast{\cal O}_\mathrm{R}(1)$ by ${\cal C}(m)$, where 
	\[ m = (d(d-1)-3)-(d+d-1-4)+1 = d^2-3d+3 = 7 , \]
	where $d=4$ is the degree of the Darboux cyclide.
	
	\begin{remark} \label{rem:sh}
		The isomorphism $(\pi|_\mathrm{R})_\ast{\cal O}_\mathrm{R}\cong{\cal C}(6)$ also implies that the global conductor ideal has a unique element of degree $6$, say $G_0$.
		The elements of degree 7 form a vector space of dimension~4, because ${\cal O}_\mathrm{R}(1)$ has four global sections (see \cite[Lemma 3.2]{GLSV}). 
		We already know a three-dimensional subspace, namely $\langle xG_0,yG_0,zG_0\rangle$.
		Let $G_1$ be an element of degree~7 that is not a multiple of $G_0$.
		Then the map $\mathrm{B}\dashrightarrow \mathrm{R}$ that is a rational inverse to $\pi|_\mathrm{R}$ can be defined as 
		\[ (x:y:z) \mapsto \left(xG_0(x,y,z):yG_0(x,y,z):zG_0(x,y,z):G_1(x,y,z)\right) = \left(x:y:z:\frac{G_1(x,y,z)}{G_0(x,y,z)}\right) . \]
	\end{remark}
	
	Let us now compute the local conductor ideals for each of the points over which $\pi|_\mathrm{R}$ is not an isomorphism.
	We will address these points as {\em special points} of the apparent contour.
	
	\begin{lemma} \label{lem:nc}
		If $q$ is a node or cusp of $\mathrm{B}_1$, and not the image of an isolated singular point, then $C_q$ is the maximal ideal.
	\end{lemma}
	
	\begin{proof}
		This is well known (see \cite[Remark~7.7]{Boehm:17} or \cite[Lemma~A.1]{GLSV}).
	\end{proof}
	
	\begin{lemma} \label{lem:cross}
		If $q\subset \mathrm{B}_1\cap \mathrm{C}_E$ is the image of two distinct points $q_1\in \mathrm{R}_1$ and $q_2\in \mathrm{C}_A$, then $C_q$ is the sum of the ideal of $\mathrm{B}_1$ and the ideal of the multiple component, i.e., it is the square of the ideal of $\mathrm{C}_E$ in the nodal case and the cube  of the ideal of $\mathrm{C}_E$ in the cuspidal case.
	\end{lemma}
	
	\begin{proof}
		This is covered by \cite[Lemma~A.1]{GLSV}. In general, if $R$ is a UFD and if $F,G$ are coprime in $R$, then the
		conductor of $R/\langle FG\rangle$ in $R/\langle F\rangle\times R/\langle G\rangle$ is $\langle F,G\rangle$. 
	\end{proof}
	
	The remaining special points are intersections of $\mathrm{B}_1$ and $\mathrm{C}_E$.
	
	A {\em covariant} is a function $F$ from the set of ideals of a regular local ring to itself such that
	any local analytic automorphism of the local ring that sends $I$ to $I'$ also sends $F(I)$ to $F(I')$.
	An example is the derivative ideal $\partial(I)$ of $I$, which is defined as the ideal generated by all partial derivatives
	of elements in $I$. 
	
	Let $I$ and $J$ be two ideals of a regular local ring, and let $\lambda\in\CC^\ast$ be a nonzero constant.
	The {\em mixed derivative ideal} $\partial(I,J,\lambda)$ is defined as the ideal generated by all elements of the form 
	$\partial(F)G+\lambda F\partial(G)$ with $F\in I$, $G\in J$, and $\partial$ being any partial derivative.
	The mixed derivative ideal is a covariant of $I$ and $J$.
	If $I=\langle F\rangle$ and $J=\langle G\rangle$, then the mixed derivative ideal is generated
	by $FG$ and all elements $\partial(F)G+\lambda F\partial(G)$, where $\partial$ is a partial derivative.
	
	Let $I,J,\lambda$ be as above, and let $K$ be another ideal. The {\em mixed jacobian ideal} $\partial(I,J,\lambda,K)$
	is defined as the ideal generated by all elements of the form
	$\frac{\partial(F,H)}{\partial(x,y)}G+\lambda F\frac{\partial(G,H)}{\partial(x,y)}$ with $F\in I$, $G\in J$, $H\in K$,
	and $x,y$ being any regular system of parameters. 
	The mixed jacobian ideal is a covariant of $I$, $J$, and $K$.
	If $I=\langle F\rangle$ and $J=\langle G\rangle$, then the mixed jacobian ideal is generated
	by the product of $FG$ and the jacobian ideal of $K$, the product $K$ and of the mixed derivative ideal $\partial(I,J,\lambda)$ and $K$,
	and all elements $\frac{\partial(F,H)}{\partial(x,y)}G+\lambda F\frac{\partial(G,H)}{\partial(x,y)}$, $H\in K$.
	
	In the following, we will use the notation $\Id(\mathrm{C}_E)$ and $\Id(\mathrm{B}_1)$ for the local ideal of the curves $\mathrm{C}_E$ and $\mathrm{B}_1$.
	
	\begin{lemma} \label{lem:touch}
		Assume that $\mathrm{C}_A$ is a nodal curve in $\mathrm{D}$.
		Assume that $q$ is the image of a point of intersection of $\mathrm{C}_A$ and $\mathrm{R}$.
		Then the local conductor is
		\[ C_q = \partial(\Id (\mathrm{B}_1),\Id (\mathrm{C}_E),-4) . \]
	\end{lemma}
	
	\begin{proof}
		This is covered by the appendix of \cite{GLSV}; see Remark~\ref{rem:newproof} for a different proof.
		% Note that the second summand $M\cdot \Id(\mathrm{C}_E)$ is not mentioned in \cite{GLSV}; this is a mistake in \cite{GLSV}. 
		
	\end{proof}
	
	\begin{lemma} \label{lem:new}
		Let $q$ be as in Lemma~\ref{lem:touch} above, but now assume that $\mathrm{C}_A$ is cuspidal.
		Let $M$ be the maximal ideal at $q$. Then the local conductor is
		\[ C_q = \partial(\Id (\mathrm{B}_1),\Id (\mathrm{C}_E),-9,M^2) . \]
	\end{lemma}
	
	\begin{proof}
		We can choose local coordinates around $p$ as in the proof of Proposition~\ref{prop:noniso}, subcase three:
		the projection map is $(x,y,z)\mapsto (x,y)$, and the surface $\mathrm{D}$ has equation $F:=z^3-(y-xz)^2=0$. The ideal of $\mathrm{B}$ is generated
		by $y^3(4x^3-27y)$, and the ideal of $\mathrm{R}$ is $\langle F,\partial_z F\rangle$. Using the computer algebra system Maple,
		we can compute the conductor and verify that it is indeed generated by the formula stated above.
		
		To prove that the formula is true in general, it is essential that the formula above is a
		covariant of $\Id(\mathrm{B}_1)$ and $\Id(\mathrm{C}_E)$.
		Let $\alpha:\CC[[x,y,z]]\to \CC[[x,y,z]]$ be an analytic coordinate change that
		transforms the local equations at $p$ to the equations above -- in particular, the subring $\CC[[x,y]]$ is mapped to itself
		because it corresponds to the projection map in both coordinate systems. Then $\alpha|_{\CC[[x,y]]}$ maps the conductor ideal
		in the coordinate system above to the conductor in the local system of coordinates at $q$. By the covariance of the formula, we may
		use the formula above also for $q$. 
	\end{proof}
	
	\begin{remark} \label{rem:newproof} \rm
		The formula for $C_q$ in Lemma~\ref{lem:touch} is covariant as well.
		Hence we can prove Lemma~\ref{lem:touch} in the same way as the proof of Lemma~\ref{lem:new} above.
		This proof is shorter and simpler than the one given in \cite{GLSV}, where an analysis by blowing ups was used.
		
		The Maple computation for computing the local conductor ideals and comparing it with the formulas
		stated above is available in \url{https://www.risc.jku.at/people/jschicho/pub/darboux/condin}
		and \url{https://www.risc.jku.at/people/jschicho/pub/darboux/condout}.
	\end{remark}
	
	The combination of all local conductor ideals to a global homogeneous conductor ideal works as follows:
	For each special point $q$, compute the maximal homogenous ideal $I_q$ such that the localization at $q$ is equal to $C_q$.
	This can be done by taking any ideal such that the localization at $q$ is equal to $C_q$ (as computed by the formulas above) and saturating it by a sufficiently high power of the maximal ideal at $q$.
	Then compute the intersection of all ideals $I_q$, where $q$ ranges over all special points.
	The final result is an ideal which has one generator $G_0$ of degree~6 and one generator $G_1$ of degree $7$.
	Indeed, the stalk of the sheaf of ideals given by the homogeneous conductor ideal at any point in $\PP^2$ is generated by the germs of dehomogenizations of
	the three polynomials $G_0$, $G_1$ and a local equation for $\mathrm{B}$.
	Hence the homogeneous conductor ideal is equal to the saturation of $\langle G_0,G_1,B\rangle$ by the irrelevant ideal $\langle x,y,z\rangle$.
	
	Once we have $G_0$ and $G_1$, we compute $\mathrm{R}$ as the image of the rational map $(x:y:z)\mapsto\left(x:y:z:\frac{G_1}{G_0}\right)$.
	If we have done everything correctly 
	(remember Remark~\ref{rem:ambi} that we need to guess, for a point in the intersection of $\mathrm{B}_1$ and $\mathrm{C}_E$, whether the point is special or whether it is a non-special image of a singular point),
	then the ideal of the image is generated by a cubic and a quartic form.

	\section{From the contour to the Surface}
	\label{sec:5}
	
	We are left with the following task: 
	given a cubic form $H_0\in\QQ[x,y,z,w]$ and a quartic form $H_1\in\QQ[x,y,z,w]$ that generate the ideal of the contour, compute the equation $F$ of the Darboux cyclide.
	Here is the information about $F$ that we can use at this point:
	
	\begin{itemize}
		\item $\partial_wF\in \langle H_0\rangle_{\QQ}$
		\item $F\in\langle xH_0,yH_0,zH_0,wH_0,H_1\rangle_{\QQ}$
	\end{itemize}
	
	The equation $F$ is only unique up to scaling, and we may make it unique by changing the first information to $\partial_wF=H_0$.
	In order to use the second information, we integrate $H_0$ formally, obtaining quartic form $H_2$ such that $\partial_wH_2=H_0$.
	Then we make an ansatz
	\[ H_2-c_1xH_0-c_2yH_0-c_3zH_0-c_4wH_0-c_5H_1 \in \QQ[x,y,z] \]
	with indeterminate symbolic coefficients $c_1,\dots,c_5$.
	This gives an inhomogeneous system of linear equations in the five indeterminates; we solve it.
	Then we set $F_0:=c_1xH_0+c_2yH_0+c_3zH_0+c_4wH_0+c_5H_1$.
	
	At this step, $F_0$ has the correct discriminant. 
	However, it is not yet the equation of a Darboux cyclide: 
	the map $\mathrm{B}\to \mathrm{R}$ we constructed is only unique up to projective transformations preserving all lines through $o$, and these do not preserve the plane at infinity.
	So, the surface $\mathrm{S}_0$ defined by $F_0=0$ is projectively equivalent to a Darboux cyclide, but its singular conic is not yet in the infinite plane.

	Let $w+ax+by+cz=0$ be the linear equation of the singular conic of $\mathrm{S}_0$. 
	Then we define $F:=F_0(x,y,z,w-ax-by-cz)$. 
	Now, the two polynomials $F$ and $F_0$ have the same discriminant with respect to $w$.
	Moreover, $F$ is the equation of a Darboux cyclide; we have solved the given problem.
	
	In Algorithm~\ref{algorithm:}, there are two issues that require an explanation. First, in lines 4 and 9, a guess is required whether a point $q$ in the apparent contour
	is the image of a special point or not. The computed result does depend on these guesses. Since we want to compute all possible Darboux cyclides that do have the given discriminant and that are in generic coordinate position, we test all possible combinations of guesses. Second, the lines 14, 15, 24, 27 are assertions that
	are checked at this point. 
	If an assertion fails, then something is wrong, and the algorithm should terminate giving an error message. The problem could be caused by a guess in line 4 or 9, so we proceed trying the next guess.
	
	\begin{proposition}
		Algorithm~\ref{algorithm:} is correct in the following sense: if the given polynomial $U$ is the equation of the apparent contour of a Darboux cyclide $\mathrm{D}$ in
		general position, then there is a guess that computes the equation of $\mathrm{D}$ up to scalings.
	\end{proposition}
	
	\begin{proof}
		Let $F'$ be the equation of a Darboux cyclide $\mathrm{D}$ in the assumption. Let $U_1$ be the factor different from $A=x^2+y^2+z^2$ 
		(as specified in the algorithm specification). Let $\mathrm{R}'$ be the zero set of $F'$ and $\partial_wF'$, i.e., the contour of $\mathrm{D}$. 
		Let $\mathrm{B}$ be the apparent contour, and let $f':\mathrm{R}'\to \mathrm{B}$ the projection map. 
		Let ${\cal C}\subset{\cal O}_\mathrm{B}$ be the sheaf of conductors of ${\cal O}_{\mathrm{B}}$ in $(f')^\ast{\cal O}_{\mathrm{R}'}$.
		Each singular point $q$ of $U_1$ is either the image of an isolated double point of $\mathrm{D}$ or not; we may assume that
		we make the correct guesses in line 4. Then, by Lemma~\ref{lem:nc}, the stalk of ${\cal C}$ at $q$ is generated by the maximal ideal at $q$.
		Similarly, each intersection point $q$ of $\mathrm{B}$ and $\mathrm{C}_E$ is either the image of a special point of $\mathrm{D}$ or not; we may assume that 
		we make the correct guess in line 9. Then, by Lemmas~\ref{lem:cross}, \ref{lem:touch}, and \ref{lem:new}, the stalk of ${\cal C}$ at $q$ is
		generated by $C_q$ as computed in line 8. After the end of the for loop in line 11, the ideal $C$ is the ideal corresponding
		to ${\cal C}$. By Remark~\ref{rem:sh}, $C$ has an element $G_0$ of degree $6$ and an element $G_1$ of degree~7 as asserted in lines 14 and 15.
		By the same remark, the rational map in line 18 maps $\mathrm{B}$ to a curve $\mathrm{R}\in\PP^3$, and the two curves $\mathrm{R},\mathrm{R}'\in\PP^3$
		are related by a projective isomorphism $\tau$ fixing the homogeneous coordinates $x,y,z$. In particular, the assertion in line 19 holds. 
		The zero set of equation $F$ computed in lines 21--23 is a quartic surface $\mathrm{S}$ that has contour $\mathrm{R}$. Such a surface, if exists,
		is necessarily unique: otherwise there would be linear pencil $F_{\lambda_1,\lambda_2}=\lambda_1 F_1+\lambda_2 F_2$ of quartics such that 
		$\partial_w F_{\lambda_1,\lambda_2}$ lies in the ideal of $\mathrm{R}$. However, this ideal only contains a single cubic form,
		therefore there exists an element $\bar{F}$ in the pencil such that $\partial_w\bar{F}=0$. This is only possible if the zero set of $\partial_w\bar{F}$ is a cone with center $o$. The projection of this cone would be
		a quartic, contradicting the fact the zero set of $U_1$ lies in this zero set. So, there is only a single quartic with
		contour $\mathrm{R}$. However, $\tau^{-1}(\mathrm{D})$ also has the 
		property that its contour is $\mathrm{R}$. This implies that $\tau$ necessarily transforms $\mathrm{S}$ to $\mathrm{D}$. 
		At this point, we see the assertion in line 24 holds; moreover, the constants $a,b,c$ are unique, because a Darboux cyclide only has one
		singular conic. The projective transformation 
		\[ \tau':\PP^3\to\PP^3, \ (x:y:z:w)\mapsto (x:y:z:w-ax-by-cz) \]
		transforms the singular conic of $\mathrm{S}$ to the singular conic of $\mathrm{D}$. The last computation in line 26 applies the transformation to $\mathrm{S}$. 
		The result is a Darboux cyclide $\mathrm{D}'$ that is projectively equivalent to $\mathrm{D}$, 
		by the equivalence $\tau'\circ\tau^{-1}$ that fixes the coordinates $x,y,z$ and that fixes the infinite plane $w=0$, because the infinite plane contains the common singular curve of both Darboux cyclides.
		But then this projective equivalence is a scaling.
	\end{proof}
	
	In order to test the method, one author took some examples of Darboux cyclides in \cite{Takeuchi:00} and applied a randomly generated Euclidean congruence transformation. 
	This author also computed the discriminant and gave it to the other authors (without revealing the surface equation). 
	The other two authors used the method in this paper, the computer algebra system Maple, and some interactive guesses of the most probable types of special points -- by Remark~\ref{rem:ambi}, there is a finite ambiguity that cannot be avoided. 
	In all cases, the secret equation could be recovered up to scaling, despite the predicted ambiguity. The Darboux cyclides that have the same apparent contour, in particular the inversion at the sphere centered at $o$, do not satisfy our assumptions on genericness of the camera position. 
	
	The computing times were less than 3 CPU seconds on a Pentium with 1.6 GHz, in any example. This is negligible in relation to the time needed for the analysis of special points and estimating their types. 
	
	\begin{example}\rm
		The input for the computation in the interactive Maple file \url{https://www.risc.jku.at/people/jschicho/pub/darboux/maplein} was computed by a random translation 
		of a random instance an equation in \cite{Takeuchi:00}. Its discriminant factors into an octic $U_1$ times $A^2$,
		which means that the absolute conic is nodal. To compute the singular points of $U_1$, we compute its discriminant with respect to $x$. 
		It factors into $D_1^3D_2^2D_3^2D_4$ with $\deg(D_1)=12$, $\deg(D_2)=4$, $\deg(D_3)=2$, and $\deg(D_4)=8$. So, we have 12 cusps, 
		and $D_1$ is the elimination ideal of the set of cusps. A Darboux cyclide cannot have 12 isolated singularities, hence all 12 cusps
		are not images of isolated singularities. The radical ideal of cusps -- which is the intersection of all maximal ideals at all cusps -- is computed 
		in one step by saturating the ideal generated by the partial derivatives of $U_1$ and $D_1$.
		
		The curve defined by $U_1$ has six nodes, four projecting to the zero set of $D_2$ and two projecting to the zero set of $D_3$. It turns out
		that the two projecting to the zero set of $D_3$ are also on $\mathrm{C}_E$. This shows that they are not images of isolated singularities,
		because a generic projection does not project isolated singularities onto $\mathrm{C}_E$. But we have to guess whether the other four nodes
		are images of isolated singularities or not.
		
		To analyze the common zeroes of $U_1$ and $A$, we compute the resultant of these two polynomials. It factors into $R_1^2D_3^2R_3$,
		with $\deg{R_1}=\deg{R_3}=4$. The zeroes of the squared factor $R_1$ correspond to tangential intersections,
		the zeroes of $R_3$ correspond to transversal intersections, and the zeroes of $D_3$ correspond to the nodes of $U_1$ that also lie
		on $\mathrm{C}_E$. These nodes are necessarily images of special points. 
		For the tangential intersections, we need to guess if they are special points projecting to tangential intersections. 
		For the transversal intersections, we have to guess if they are images of pinch points. 
		In total, there are 8 possible guesses. 
		
		We compute the conductor ideal for all 8 cases. The assertion in line 14 fails for 6 cases. The assertion in line 15 fails for one of the two
		remaining cases. So, we have only a single case still to consider. In this case, we have no isolated singularities and two special points,
		namely those projecting to nodes of $U_1$ that also lie on $A$.
		The result is the Darboux cyclide that was used for constructing the input, up to scaling.
		The output is shown in \url{https://www.risc.jku.at/people/jschicho/pub/darboux/mapleout}.
	\end{example}
	
	\vfill
	
	\renewcommand{\thealgorithm}{DarbouxFromApparentContour} % Removes algorithm numbering
	\begin{algorithm}[H]
		\caption{}\label{algorithm:}
		\begin{algorithmic}[1]
			\Require A polynomial $U=U_1(x^2+y^2+z^2)^2$ or $U=U_1(x^2+y^2+z^2)^3$ in $\QQ[x,y,z]$,
			the equation of the apparent contour $\mathrm{B}\subset {\PP^2}$.
			\Ensure The equation $F\in\QQ[x,y,z,w]$ of a Darboux cyclide $\mathrm{D}\subset {\PP^3}$, such that its
			discriminant is equal to $U$.
			\Statex
			\State {\bf Compute} the singular points of the zero set of $U_1$ and common zeroes of $U_1$ and $x^2+y^2+z^2$.
			\State {\bf Initialize} the conductor ideal $C:=\langle 1\rangle_{\QQ[x,y,z]}$.
			\For{all singular points $q$ of the zero set of $U_1$}
			\If{$q$ does not look like the image of an isolated double point}
			\State $C:=C\ \cap$ maximal ideal at $q$
			\EndIf
			\EndFor
			\For{all common zeroes $q$ of $U_1$ and $x^2+y^2+z^2$}
			\If{$q$ does not look like the image of a special point}
			\State Compute $C_q$ by the formulas in Lemmas~\ref{lem:cross}, \ref{lem:touch}, and \ref{lem:new}
			\State $C:=C\ \cap C_q$ 
			\EndIf
			\EndFor
			\State {\bf Assert} that $C$ has a unique element $G_0$ in degree $6$ up to scalar multiplication
			\State {\bf Assert} that $C$ has a unique element $G_1$ in degree $7$ up to scalar multiplication
			\State $\ \ $ and up to addition of multiples of $G_0$
			\State {\bf Compute} the image $\mathrm{R}$ of $\mathrm{B}$ under the rational map
			\State $\ \ $ $(x:y:z)\mapsto (xG_0:yG_0:zG_0:G_1)$
			\State {\bf Assert} that the ideal of $\mathrm{R}$ is generated by a cubic $H_0$ and a quartic $H_1$
			\State {\bf Compute} $H_2:=\int_w H_0$ which is homogeneous
			\State {\bf Find} $c_0,\dots,c_5$ such that
			\State $\ \ $ $H_2-c_1xH_0-c_2yH_0-c_3zH_0-c_4wH_0-c_5H_1 \in \QQ[x,y,z]$
			\State $F_0:=c_1xH_0+c_2yH_0+c_3zH_0+c_4wH_0+c_5H_1$
			\State {\bf Assert} that there exist $a,b,c\in\QQ$ such that $w+ax+by+cz$
			\State $\ \ $is the linear equation of a conic which is singular in $F_0$.
			\State $F:=F_0(x,y,z,w-ax-by-cz)$.
			\State \Return $F$.
		\end{algorithmic}
	\end{algorithm}

	\begin{figure}[H]\label{fig:example1}
		\includegraphics[width=0.4\textwidth]{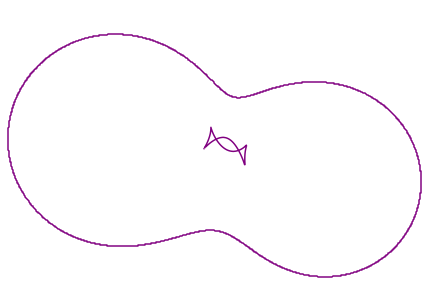} \qquad \qquad
		\includegraphics[width=0.4\textwidth]{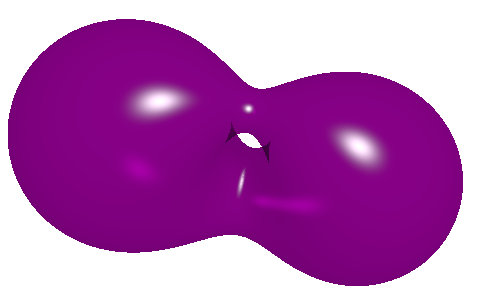} \\
		%\end{figure}
		%\begin{figure}\label{fig:example2}
		\includegraphics[width=0.4\textwidth]{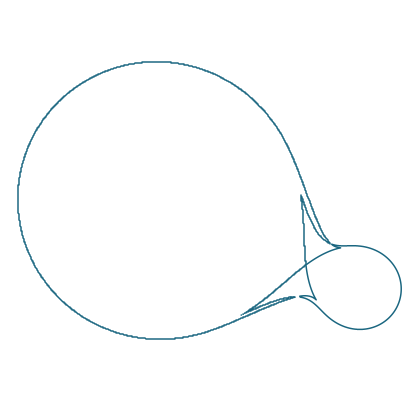} \qquad \qquad
		\includegraphics[width=0.4\textwidth]{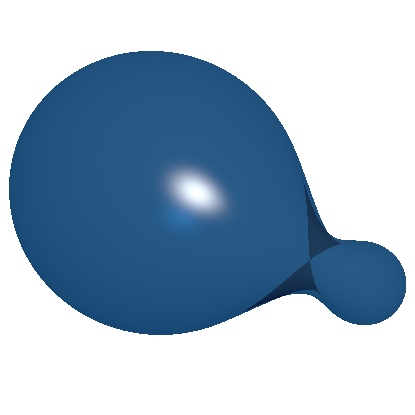} \\
		%\end{figure}
		%\begin{figure}\label{fig:example3}
		\includegraphics[width=0.4\textwidth]{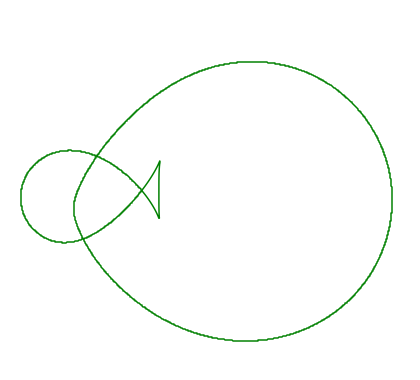}\qquad \qquad
		\includegraphics[width=0.4\textwidth]{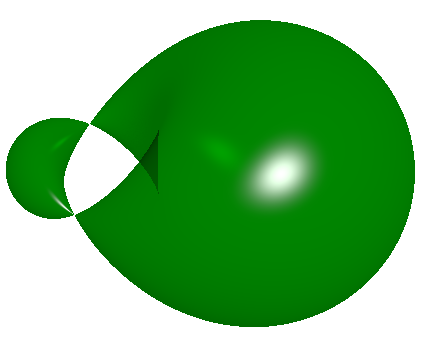} 
		\caption{Some Darboux cyclides with $\mathrm{C}_A$ being a nodal curve, together with the real part of their apparent contour, a curve of degree~8. 
			The complex part of the apparent contour is not visible: it is the elliptic absolute conic $\mathrm{C}_E$. 
			The equation of $\mathrm{C}_E$ can be assumed to be known because of the assumption that we have a calibrated camera.}
	\end{figure}
	\begin{figure}[H]\label{fig:example4}
		\includegraphics[width=0.4\textwidth]{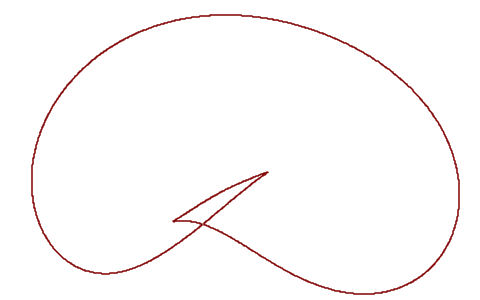} \qquad \qquad
		\includegraphics[width=0.4\textwidth]{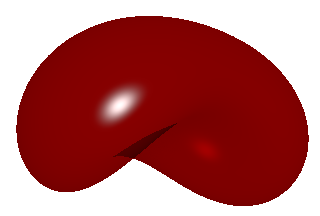} \\
		%\end{figure}
		%\begin{figure}\label{fig:example5}
		\includegraphics[width=0.4\textwidth]{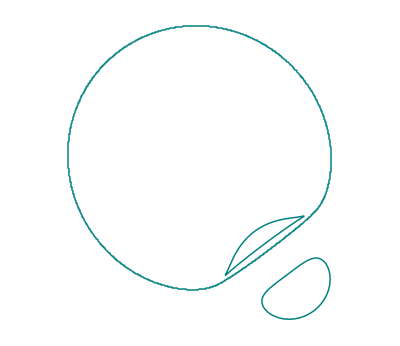} \qquad \qquad
		\includegraphics[width=0.4\textwidth]{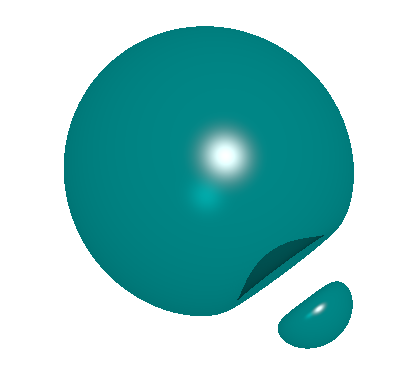} \\
		%\end{figure}
		%\begin{figure}\label{fig:example6}
		\includegraphics[width=0.35\textwidth]{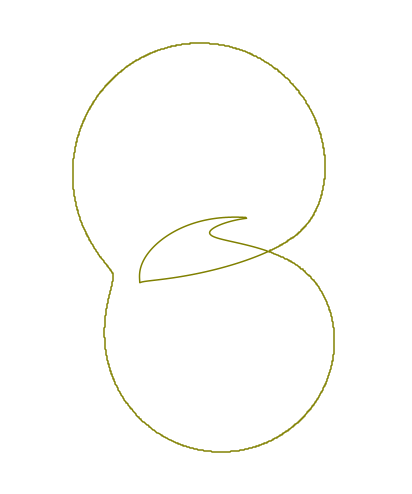} \qquad \qquad
		\includegraphics[width=0.4\textwidth]{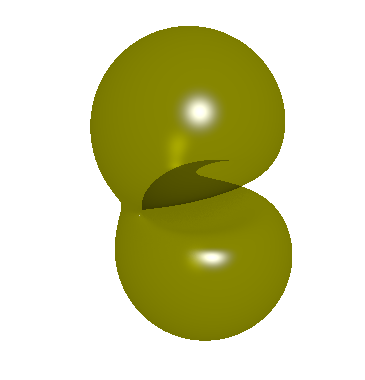} 
		\caption{Some Darboux cyclides with $\mathrm{C}_A$ being a cuspidal curve, together with the real part of their apparent contour, a curve of degree~6.}
	\end{figure}

	\section*{Acknowledgements}
	This work is part of GRAPES project that has received funding from the European Union’s Horizon 2020 research and innovation programme under the Marie Skłodowska-Curie grant agreement No 860843.
	We appreciate the discussion with Niels Lubbes for explaining the possible isolated singularities of a Darboux cyclide. 
	We thank Rimvydas Krasauskas for the suggestion on the quick online visualization \cite{CindyJS:04}, in order to choose good shapes for the pictures, and also for pointing out the other possible solution of our algorithm by applying inversion in the camera position. The pictures are produced using \cite{PovRay:04}.
	
	\bibliographystyle{unsrt}
	\small
	\bibliography{darboux}
	
\end{document}